\documentclass[11pt,a4paper]{article}
\pdfoutput=1

\usepackage[numbers,sort&compress]{natbib}
\usepackage{sjtu_iic_arxiv}
\usepackage{mathtools}
\usepackage{booktabs}
\usepackage{multirow}
\usepackage{algorithm}
\usepackage{algorithmic}
\usepackage{float}
\usepackage{placeins}
\usepackage{amssymb}

\numberwithin{equation}{section}

\newcommand{\methodname}{InterOPT}
\newcommand{\methodfullname}{Interactive Optimization}
\newcommand{\benchmarkname}{OR-Clarify}
\newcommand{\readytoken}{\texttt{READY\_TO\_MODEL}}

\hypersetup{
  pdftitle={Ask Before You Optimize: Dynamic Pre-Formulation Clarification for Interactive Optimization},
  pdfauthor={Sihan Ge, Yichen Lin, Chenyu Zhou, Jianghao Lin, Tao Yao, Dongdong Ge}
}

\githublink{https://github.com/AIOR-Research/InterOpt}
\huggingfacelink{https://huggingface.co/datasets/AIOR-Research/OR-Clarify}

\setheadertext{%
  \raisebox{-0.55cm}{\includegraphics[height=1.3cm]{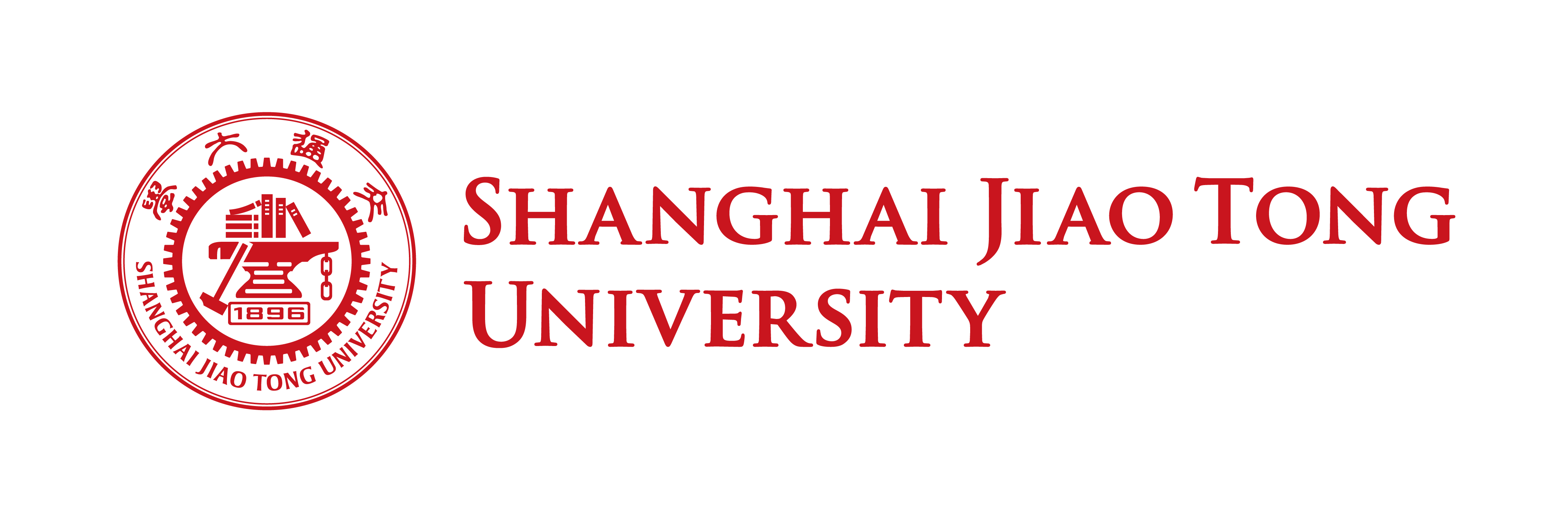}}%
  \hspace{0.15cm}%
  \raisebox{-0.29cm}{\includegraphics[height=0.72cm]{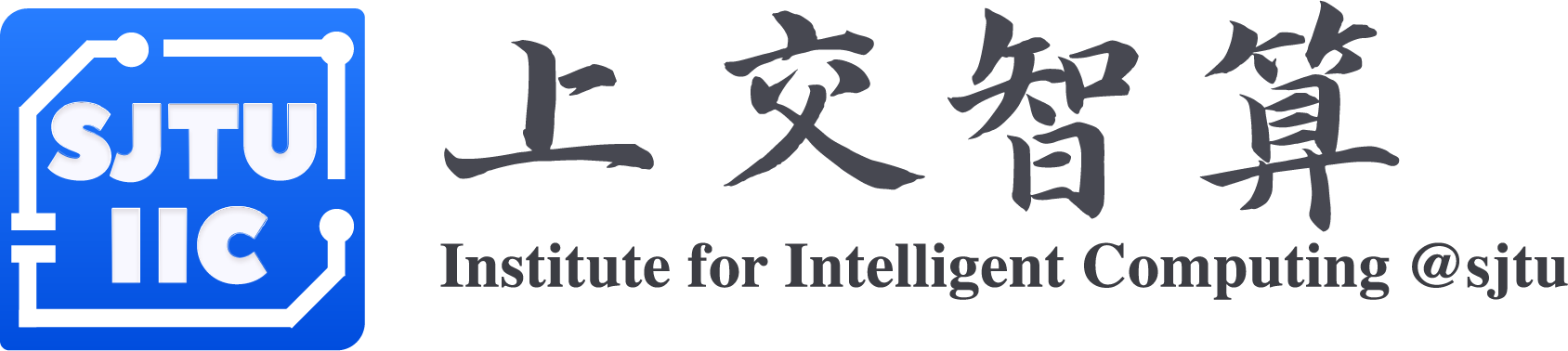}}%
}
\setheadertitle{Ask Before You Optimize: Dynamic Pre-Formulation Clarification for Interactive Optimization}

\title{%
  Ask Before You Optimize: Dynamic Pre-Formulation Clarification for Interactive Optimization
}

\author{%
  Sihan Ge$^{1,*}$,
  Yichen Lin$^{2,*}$, 
  Chenyu Zhou$^{2,*}$, 
  Jianghao Lin$^{2,\dagger}$,
  Tao Yao$^{2,\dagger}$,
  Dongdong Ge$^{2}$ \\
  {\sjtuiicAffilFont
  $^{1}$ Cardinal Operations\\
  $^{2}$ Shanghai Jiao Tong University, Shanghai, China\\
  \texttt{enoch-n@outlook.com} \quad
  \texttt{\{linyichen, chenyuzhou, linjianghao, taoyao, ddge\}@sjtu.edu.cn}\\
  $^*$ Equal contribution \quad
  $^\dagger$ Corresponding authors
  }
}

\begin{document}

\begin{abstract}
Large language models (LLMs) are increasingly used to formulate optimization models from natural-language problem descriptions, yet realistic operations research (OR) requests are often incomplete: missing objectives, constraints, or business rules can change the resulting mathematical program. Existing evaluations largely assume a complete specification and therefore overlook whether an agent knows when clarification is needed before modeling.
We introduce \benchmarkname{}, a benchmark for pre-formulation clarification. Each task presents a partial public problem description, withholds structured hidden slots, and evaluates agents through bounded interaction with a simulated user. The benchmark supports both open-ended and choice-based clarification, and measures slot recovery, stopping behavior, silent assumptions, and interaction cost.
We further propose \methodfullname{} (\methodname{}), a two-stage framework that identifies unresolved formulation-critical gaps and uses them to guide whether to ask the next question or to stop. In our choice-based experiments, \methodname{} substantially outperforms all baselines in exact slot recovery; in the open-ended setting, it remains competitive with strong prior methods. Together, OR-Clarify and InterOPT reframe OR assistance as a selective completeness decision: clarify when needed, stop when ready, and quantify what remains missing.
\end{abstract}

\maketitle

\section{Introduction}

The structure of an operations research (OR) model—including its objective function, constraints, decision variables, and feasible region—is determined by the business problem specification.
With large language models (LLMs), non-experts can now describe problems in natural language and receive a candidate mathematical formulation \citep{ramamonjison2022nl4opt,jiang2024llmopt}, lowering the barrier to OR modeling.
However, real-world business requests are rarely delivered as complete textbook-like specifications.
A user may describe capacities and demands without specifying the objective, or state a routing rule without clarifying whether the time windows are hard or soft.
These omissions are not superficial: they can fundamentally alter the mathematical structure of the problem.

The central failure mode in this setting is \textbf{\emph{premature formulation}}.
Most evaluations of LLM-based optimization agents assume a sufficient specification and measure whether the agent can solve or express a given model \citep{ahmaditeshnizi2024optimus,huang2024orlm}, missing the earlier question of whether the available information supports a meaningful formulation.
We find that strong LLM agents often declare readiness while core business facts remain unclarified, or silently fill missing facts with unsupported defaults.
For example, in a vehicle-routing context, if a user omits whether a courier must return to the origin, an agent that silently assumes a closed tour changes the constraint structure without asking for confirmation.


Recent interactive OR systems have begun to incorporate user interaction. ORPilot \citep{xie2026orpilot} uses interviews within an end-to-end modeling pipeline, while Drossman et al. \citep{drossman2026conversation} study conversational optimization through iterative solution refinement toward stakeholder utility. Taken together, these studies demonstrate the value of interaction in OR problem solving. However, interaction is embedded within a broader interactive process, and its effectiveness is not isolated as an explicit evaluation target. In particular, it remains unclear whether an agent actually recovers formulation-critical missing requirements and whether it knows when enough information has been obtained to proceed. These capabilities matter because unresolved requirements can change the resulting optimization formulation. We therefore study pre-formulation clarification as a distinct research problem and develop both a dedicated evaluation framework and clarification methods tailored to this problem.


This paper studies \textbf{\emph{pre-formulation clarification}} as a standalone OR task.
We define a fact as \textbf{\emph{formulation-critical}} if its value can change the structure of the resulting optimization formulation. The task is to determine whether the current public specification is model-ready and, when it is not, to recover the missing formulation-critical facts with as little interaction as possible. 
To address this challenge, we propose Interactive Optimization (\textbf{InterOPT}), a two-stage framework that separates gap diagnosis from interaction control. 
The first stage, \emph{Dynamic Gap Search}, identifies formulation-critical gaps and maintains a cross-turn record of those that remain unresolved by the public transcript. 
The second stage, \emph{Gap-Guided Action Search}, uses the currently unresolved gaps to decide whether to ask a targeted question or stop. By separating persistent gap tracking from action selection, InterOPT aims to improve specification recovery without unnecessary interaction.


To systematically study this task, we build \textbf{\emph{\benchmarkname{}}}, an evaluation framework and benchmark for pre-formulation clarification in OR. 
Each instance pairs an incomplete public brief with private, source-supported formulation-critical facts,
fact-bounded simulated-user responses, and slot-level recovery rubrics, enabling controlled evaluation under both free-form and choice protocols.
Its construction pipeline further converts fully specified optimization tasks into clarification instances by withholding formulation-critical facts and generating the corresponding interaction and evaluation artifacts. 
To our knowledge, \benchmarkname{} is the first OR-specific framework to jointly evaluate formulation-gap recovery, readiness decisions, and interaction cost before LLM-driven autoformulation.

\begin{figure*}[!t]
\centering
\includegraphics[width=1\textwidth]{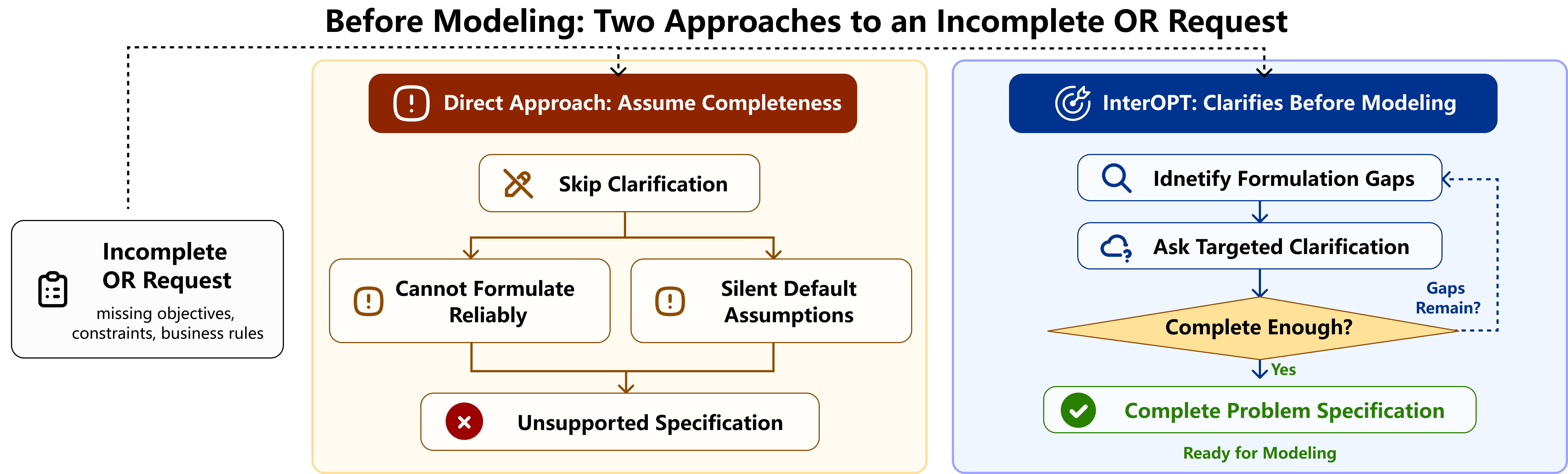}
\caption{Motivation for pre-formulation clarification.
}
\label{fig:paper_overview}
\end{figure*}

In summary, this paper makes three contributions.

\textbf{(1) The pre-formulation clarification task and the \methodname{} framework.}
We formulate pre-formulation clarification as the joint problem of assessing whether a public specification is model-ready and recovering missing formulation-critical facts while minimizing unnecessary interaction. We then propose \methodname{}, a two-stage framework in which Dynamic Gap Search identifies and tracks unresolved formulation gaps across turns and Gap-Guided Action Search uses this state to decide what to ask and when to stop.

\textbf{(2) The \benchmarkname{} construction and evaluation framework.}
\benchmarkname{} converts fully specified optimization tasks into controlled clarification instances by withholding formulation-critical facts and generating the corresponding public briefs, fact-bounded simulated-user responses, and slot-level evaluation artifacts. It supports controlled evaluation under both free-form and Choice interaction protocols.

\textbf{(3) A systematic empirical study of pre-formulation clarification.}
We compare \methodname{} with strong LLM-based baselines and analyze exact requirement recovery, readiness decisions, silent assumptions, and interaction cost. The results demonstrate strong recovery gains in the Choice setting and competitive performance in the free-form setting, while the ablations and behavioral diagnostics reveal the roles of gap tracking, question selection, and stopping.
\section{Related Work}

\subsection{LLMs for Optimization Modeling}

Recent work studies whether LLMs can translate natural-language problem descriptions into optimization models and solver-ready code.
Benchmarks and systems such as NL4Opt, OptiMUS, LLMOPT, and ORLM evaluate model generation, solver integration, and domain adaptation \citep{ramamonjison2022nl4opt,ahmaditeshnizi2024optimus,jiang2024llmopt,huang2024orlm}, and production-oriented systems such as ORPilot organize modeling into structured pipelines spanning interview, data collection, code generation, and execution \citep{xie2026orpilot}.
Most prior benchmarks assume complete specifications. ORPilot supports clarification through interviews, but does not directly evaluate hidden-slot recovery or readiness.

Training-oriented work further improves OR-specific reasoning through process supervision. StepORLM introduces a self-evolving framework that combines solver-based outcome verification with generative process supervision \cite{zhou2025steporlmselfevolvingframeworkgenerative}. Recent agentic methods also explore coordination and structured workflow construction: Agora-Opt combines decentralized debate with a reusable memory bank \cite{lin2026soliloquyagoramemoryenhancedllm}, while LEAN-LLM-OPT constructs structured multi-agent workflows for large-scale optimization auto-formulation \cite{liang2026largescaleoptimizationmodelautoformulation}. On the evaluation side, OR-Space broadens OR-agent benchmarking beyond one-shot formulation by evaluating Build, Revise, and Explain tasks in persistent multi-artifact workspaces \cite{zhou2026orspacefulllifecycleworkspacebenchmark}. These developments expand OR-oriented training, coordination, workflow design, and lifecycle evaluation. Pre-formulation clarification nevertheless remains underexplored as a standalone target: whether the current specification is model-ready, which formulation-critical facts are still missing, and when the agent should stop asking.

\subsection{Clarification, Elicitation, and Abstention}

Clarifying-question research studies how agents ask questions when user intent is underspecified, with datasets and methods for conversational retrieval \citep{aliannejadi2019qulac,aliannejadi2020clariq}, ambiguity resolution in open-domain QA \citep{min2020ambigqa}, and selective clarification \citep{kuhn2022clam}.
Conversational machine reading makes missing rule conditions explicit and permits follow-up questions before a decision \citep{saeidi2018sharc}, while CAmbigNQ represents alternative interpretations through a clarification question with user-selectable options \citep{lee2023cambignq}. Although these settings establish both free-form and option-based clarification, they primarily target search intents, answers, or rule conditions, leaving requirements that modify an optimization formulation unaddressed.

Preference-elicitation work asks informative questions to improve downstream decisions \citep{tamkin2023gate}, teaches models to ask better clarifying questions \citep{andukuri2024stargate}, and optimizes multi-turn trajectories \citep{dou2025togate,zhang2025futureturns}; more recent work treats clarification as a decision about when to ask, what to ask, and when to stop. \citet{zhang2025clarify} introduce IntentSim, which estimates the value of clarification from the entropy over simulated user intents; AskBench evaluates missing-intent and false-premise settings with an interactive judge and simulated user \citep{zhao2026askbench}; and SAGE-Agent selects questions from structured uncertainty and introduces ClarifyBench for multi-turn tool disambiguation \citep{suri2026structured}. In the optimization context, \citet{drossman2026conversation} showed that conversational interaction improves solution quality over one-shot submission.
These methods provide close points of comparison because they treat questioning as active information gathering. Still, their targets are general user intent or preferences, whereas a question here is correct only if it recovers a business fact that can change the optimization formulation.

Adjacent interactive benchmarks increasingly evaluate whether dialogue reaches a correct formal or environmental state. \mbox{$\tau$-bench} couples tool-using agents with simulated users and scores the final database state and cross-run reliability \citep{yao2025taubench}. CLARITY constructs single- and multi-turn ambiguity cases for NL2SQL and evaluates localization and resolution of schema-level ambiguity \citep{sarwar2026clarity}. In contrast, \benchmarkname{} makes formulation-critical business requirements the private targets and evaluates their complete recovery, silent assumptions, and readiness before an optimization model is built.

A separate line of work evaluates whether LLMs can abstain when information is insufficient \citep{kirichenko2025abstentionbench,kadavath2022language}.
Generic abstention asks, "Can I answer this question?" In contrast, we ask, "Is the current information sufficient to build the correct model?"---which requires identifying the missing fact that changes the formulation, obtaining it, and stopping only after core slots are recovered.
This makes readiness a joint problem of uncertainty detection, assumption localization, and interaction control.
\section{Problem Formulation}
\label{sec:problem_formulation}

We study \textbf{\emph{pre-formulation clarification}}: the task of
(1) deciding whether a business request contains enough information to determine a meaningful optimization formulation, and
(2) recovering formulation-critical information that remains unresolved before formulation begins.

Let $b$ denote the initial public business brief, and let
$\tau_t=\big((a_s,u_s)\big)_{s<t}$
denote the public interaction transcript before turn $t$, where $s$ indexes preceding interaction turns, $a_s$ is the agent's public clarification action, and $u_s$ is the user's response.
Let $P_t$ denote the public problem statement induced by $b$ and $\tau_t$.

Let $\mathcal{B}(P_t)$ be the set of plausible business completions consistent with $P_t$.
For $c\in\mathcal{B}(P_t)$, let $\phi(P_t,c)$ denote the formulation structure induced by completing $P_t$ with $c$, including the objective, constraints, decision variables, and other formulation-level structures.
We say that $P_t$ is \emph{formulation-complete} if and only if all plausible completions induce the same formulation structure:
\begin{equation}
\left|\{\phi(P_t,c): c\in\mathcal{B}(P_t)\}\right|=1.
\end{equation}
Conversely, $P_t$ is \emph{formulation-incomplete} if plausible completions can induce more than one formulation structure:
\begin{equation}
\left|\{\phi(P_t,c): c\in\mathcal{B}(P_t)\}\right|>1.
\end{equation}
Formulation incompleteness therefore refers specifically to unresolved business conditions whose alternative resolutions can change the induced optimization formulation.

At each turn, the agent either asks a clarification question or judges that the current information is sufficient for formulation.
If the agent asks, the user's response is appended to the public transcript and induces an updated problem statement $P_{t+1}$.
The interaction terminates when the agent emits \readytoken{} or reaches the maximum turn limit $T_{\max}$.

Pre-formulation clarification requires the agent to determine what information should be requested next to resolve formulation-relevant uncertainty. It must also judge when the current public information is sufficient to proceed with formulation. A successful clarification policy should recover the information needed to reach a formulation-complete state while avoiding unnecessary interaction and premature readiness.
\section{OR-Clarify Benchmark and Evaluation Framework}
\label{sec:benchmark}

To evaluate \textbf{\emph{pre-formulation clarification}} under controlled conditions, \benchmarkname{} turns formulation-relevant uncertainty into explicit, benchmark-side targets.
Each case $i$ represents the setting as a public--private tuple
$(b_i,\mathcal{F}_i,\mathcal{H}_i)$, where $b_i$ is the public business brief, $\mathcal{F}_i$ is the source-grounded fact set withheld from the agent, and $\mathcal{H}_i=\{h_{i,j}\}_{j=1}^{m_i}$ is the set of hidden slots that are grounded in $\mathcal{F}_i$ but absent from $b_i$.
Throughout, $i$ indexes cases, $m_i$ denotes the number of hidden slots in case $i$, and $j\in\{1,\ldots,m_i\}$ indexes those slots.

These hidden slots serve as finite, benchmark-side annotations of unresolved conditions that can change the resulting formulation defined in Section~\ref{sec:problem_formulation}.
An agent's clarification ability is then measured by whether it actively recovers the underlying requirements through interaction and declares readiness only after the relevant information has been established.

\subsection{Benchmark Construction Framework}

Benchmark construction begins with a complete, source-grounded record of an intended OR task.
We first decompose each record into individual facts describing the business setting, numerical inputs, objective, operational constraints, and modeling assumptions. Each fact expresses a single requirement.

We then determine which facts may be withheld from the initial public brief.
The business setting and numerical inputs remain visible, while facts about the objective, constraints, or assumptions are eligible for masking.
A candidate fact is excluded from masking when its mathematical meaning is already determined by the visible information. For instance, in a multi-period production-planning problem, 15,000 available production hours already implies a capacity upper bound, while a demand of 1,000 units does not determine whether it must be met exactly or may be backlogged. The latter leaves a demand-satisfaction rule that may require clarification.
This screening step retains only information gaps that can meaningfully affect the formulation.

Among the eligible facts, a deterministic pseudorandom procedure with a fixed, case-specific seed masks roughly half.
This masking rate is fixed before method evaluation, yielding briefs that are partially specified yet still interpretable.
Eligible facts that are not selected remain visible, and all visible facts are compiled into a self-contained public brief $b_i$ shown to the agent.
The selection is then frozen so that every evaluated method starts from the same information and faces the same missing requirements.

Each masked fact becomes a single \emph{hidden slot}, representing one missing requirement against which the agent's clarification behavior is evaluated.
Each slot is linked one-to-one to its underlying fact and includes supporting evidence, a simulated-user answer restricted to that fact, examples of acceptable questions, a semantic recovery rule, and a severity label.
These annotations allow the judge to recognize semantically equivalent successful questions without requiring an exact wording match.
Together, the hidden slots in case $i$ form its evaluation target set $\mathcal{H}_i$.

After masking, each hidden slot is assigned a severity label in
$\{\mathrm{P0},\mathrm{P1},\mathrm{P2}\}$.
A $\mathrm{P0}$ slot denotes a blocking condition whose omission can change the problem itself, such as whether a route is open or closed.
A $\mathrm{P1}$ slot denotes a substantive modeling condition whose omission can leave the formulation incomplete or materially incorrect, such as whether unmet demand is penalized.
A $\mathrm{P2}$ slot denotes a secondary boundary or interpretive condition that affects modeling fidelity and is less central to the core evaluation, such as whether vehicles may be scheduled across day boundaries.
Because severity is assigned only after masking, severity labels do not influence the masking procedure.

The same construction procedure, consisting of decomposition, screening, masking, and annotation, can be applied to additional complete, source-grounded OR task records.
Applying it to our current source collection yields \benchmarkname{}, which comprises 100 clarification cases and 178 hidden slots, with 1--5 slots per case (mean 1.78): 75 $\mathrm{P0}$, 83 $\mathrm{P1}$, and 20 $\mathrm{P2}$ slots.
Human auditing verifies slot boundaries, severity labels, answer support, and rubric consistency.

\subsection{Controlled Information Boundary}

\benchmarkname{} enforces a strict information boundary.
During interaction, the tested agent sees only the public brief $b_i$ and the public transcript.
The simulated user has access to the private case facts $\mathcal{F}_i$ but answers only the current question and never volunteers unasked hidden facts.

Under the Choice setting, the simulated user selects among the available options based on the private case facts.
In MC-D-based variants, it selects option D with a short free-form correction whenever none of A--C is supported by the private facts.
Its rationale and match label are used only for auditing and are withheld from the tested agent.

Once the interaction ends, the judge receives the frozen hidden-slot annotations and the public transcript.
A slot receives exact-recovery credit only if, before \readytoken{}, an agent question or an explicit assumption check semantically identifies that requirement.
Facts volunteered without being requested, assumptions introduced only in the final model, vague catch-all questions, and partial matches receive no exact credit.
For every slot, the judge records the supporting transcript location along with a \texttt{yes}, \texttt{partial}, or \texttt{no} label.

A separate protocol detector checks whether the public actions follow the required interaction format and supplies no recovery information.
This separation keeps the hidden slots and evaluation rubrics strictly outside the tested agent's information boundary.

\benchmarkname{} supports two interaction settings: an open/free-form setting, where the agent asks natural-language clarification questions, and a Choice setting, where the agent poses questions with candidate options.

\begin{figure*}[!t]
\centering
\includegraphics[width=1\textwidth]{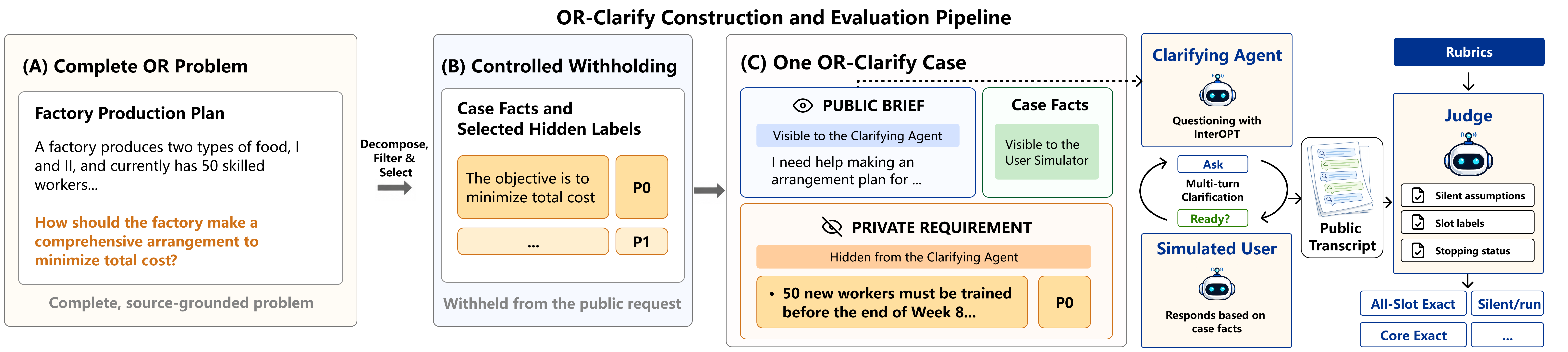}
\caption{Evaluation workflow and information flow. Hidden slots never enter the evaluated agent. The simulator supplies controlled answers, the passive monitor records protocol behavior, and the post-hoc judge scores the completed public transcript against frozen hidden-slot rubrics.}
\label{fig:evaluation_workflow}
\end{figure*}

Figure~\ref{fig:evaluation_workflow} summarizes the benchmark construction, the controlled information boundary, the interaction workflow, and the post-hoc evaluation process.

\subsection{Evaluation Protocol and Metrics}

The evaluation protocol runs $N$ cases with $K$ repeated runs per case and at most $T_{\max}$ turns per run.
All metrics are computed per run, averaged within each case, and then averaged across cases.

We define the core hidden-slot set as: 
\begin{equation}
\mathcal{H}_i^{\mathrm{core}}
=
\{h_{i,j}\in\mathcal{H}_i:
h_{i,j}\text{ has severity }\mathrm{P0}\text{ or }\mathrm{P1}\}.
\end{equation}
Cases that contain no $\mathrm{P0}/\mathrm{P1}$ slots are excluded from core-based evaluations but remain part of the full-slot analysis.
Let $\mathcal{I}_{\mathrm{core}}=\{i:|\mathcal{H}_i^{\mathrm{core}}|>0\}$ denote the set of cases containing at least one core hidden slot.

For run $k$ of case $i$, let $z_{i,k,j}=1$ if hidden slot $h_{i,j}$ is exactly recovered in the final public transcript and $0$ otherwise.
Partial recovery does not count as exact recovery.

Core Exact is the primary core-completeness metric.
A run succeeds under this metric only when all $\mathrm{P0}/\mathrm{P1}$ hidden slots in a core-eligible case are exactly recovered:
\begin{equation}
\mathrm{CoreExact}
=
\frac{1}{|\mathcal{I}_{\mathrm{core}}|}
\sum_{i\in\mathcal{I}_{\mathrm{core}}}
\frac{1}{K}
\sum_{k=1}^{K}
\mathbf{1}\!\left[
\forall h_{i,j}\in\mathcal{H}_i^{\mathrm{core}},\ 
z_{i,k,j}=1
\right].
\end{equation}

All-Slot Exact applies the same criterion to every hidden slot, including $\mathrm{P2}$:
\begin{equation}
\mathrm{AllExact}
=
\frac{1}{N}
\sum_{i=1}^{N}
\frac{1}{K}
\sum_{k=1}^{K}
\mathbf{1}\!\left[
\forall h_{i,j}\in\mathcal{H}_i,\ 
z_{i,k,j}=1
\right].
\end{equation}

We report several diagnostics alongside exact recovery.
A silent assumption is recorded when a hidden requirement has not been confirmed, yet the agent later treats one particular value as established in a clarification turn, its readiness summary, or its final answer.
For example, stating that every route returns to the depot without first checking whether routes are open or closed constitutes a silent assumption. Asking that question or explicitly listing the issue as unresolved does not.

The judge also labels each run's stopping behavior as premature, appropriate, over-questioning, or no-stop.
Interaction burden is reported using Avg Turns and Avg Q. Avg Q counts atomic clarification questions and may exceed Avg Turns when a single turn contains several questions.
Together, these metrics assess requirement recovery, readiness behavior, silent assumptions, and interaction efficiency under a single controlled information boundary.
\section{Method}

\subsection{Overview}

We introduce \methodfullname{} (\methodname{}), a two-stage framework for formulation-gap-guided clarification prior to optimization.
Its two stages separate two decisions that are easy to conflate in LLM-based clarification: diagnosing what is still missing, and choosing what to ask next.

The decomposition is motivated by a monitoring-control view of problem solving.
In this view, a reasoner not only performs task actions but also monitors what is known, what is still uncertain, and whether the current state is sufficient for the next step \citep{nelson1990metamemory}.
Control processes then use this monitored state to allocate further search or to stop. We use this distinction as a design lens. 

In pre-formulation OR clarification, the monitoring problem is to identify missing facts that could change the optimization formulation if resolved differently. A request can look model-ready while omitting an objective convention, a feasibility rule, a decision boundary, or a hard-versus-soft policy.

Directly generating the next plausible question can miss the deeper issue of whether the agent has diagnosed the right gap. \methodname{} turns monitoring into a concrete two-stage loop. Stage 1, Dynamic Gap Search, continuously identifies formulation gaps and maintains them in a persistent ledger. Stage 2, Gap-Guided Action Search, uses the currently open entries, when available, to generate three gap-bound candidate questions, from which a selector chooses one to pose. The Stage 2 action generator $\Psi$, instantiated by the tested model, separately decides whether to ask or declare readiness.

When the open set $\mathcal{O}_t$ is non-empty, every generated candidate question is anchored to a specific open gap, keeping question generation aligned with the diagnostic memory from Stage 1. 
The stopping decision, however, remains with the model and is not determined by the ledger. 
At each turn, $\Psi$ emits either \textsc{Ask} or \readytoken{}. 
It is instructed to favor \textsc{Ask} when an unresolved assumption could change the objective, decision scope, constraints, entities, or time structure, and to allow \readytoken{} when the remaining uncertainty concerns only notation, raw-data collection, or downstream solver details.

\begin{figure}[!t]
\centering

\includegraphics[width=\textwidth]{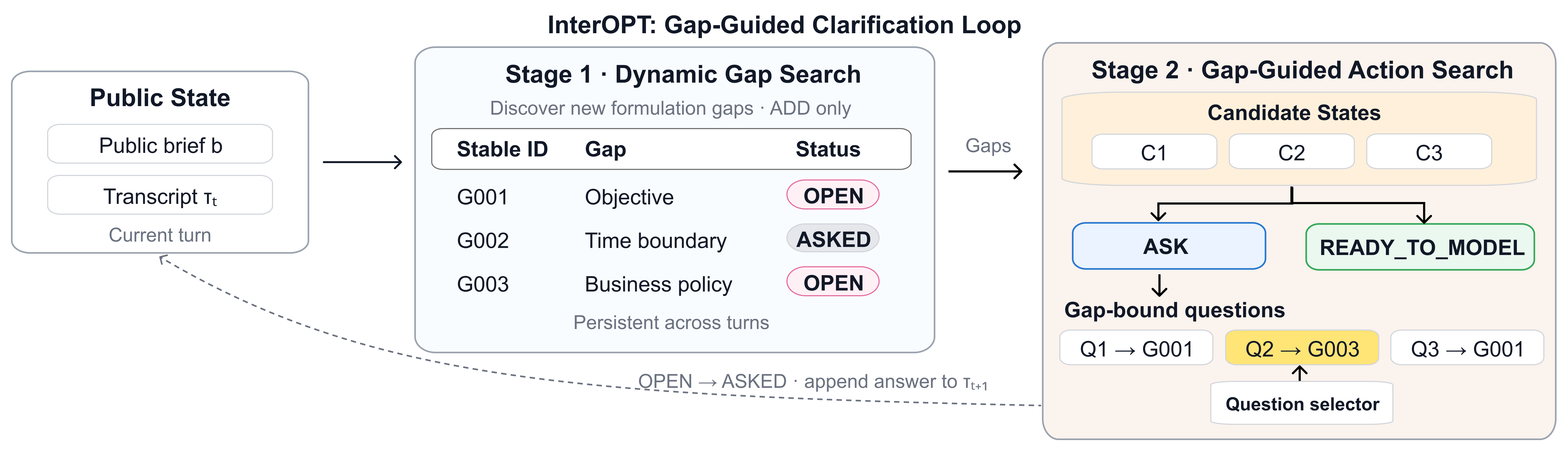}
\caption{\methodname{} gap-guided clarification loop. Stage 1 maintains a persistent public-evidence gap ledger, and Stage 2 selects a gap-bound question or \texttt{READY\_TO\_MODEL}.}
\label{fig:interopt_pipeline}
\end{figure}

This separation yields two properties. First, the runner never blocks \readytoken{} merely because open entries remain: it terminates the interaction and records the event for diagnosis. Second, an empty open set does not force readiness, since Stage 1 may have missed a gap. Neither stage has access to the benchmark's hidden slots or their P0--P2 labels. The ledger therefore provides persistent memory and question--gap alignment for the uncertainties it captures, but it does not guarantee complete discovery or safe stopping.

A \emph{formulation gap} is an unresolved business condition whose different plausible values would lead to different optimization formulations.
Stage 1 discovers and tracks these gaps; Stage 2 acts on them.
The method operates solely on the public brief and the dialogue history; it does not access benchmark-side hidden slots, private simulator facts, or evaluation rubrics.

\subsection{Stage 1: Dynamic Gap Search}

Before choosing the action at turn $t$, let $\mathcal{M}_t$ denote the agent's internal gap memory.
Each entry in $\mathcal{M}_t$ represents a formulation-critical business condition whose value is not fully determined by the public brief and transcript.
The memory is derived solely from public evidence; it does not contain benchmark hidden slots or simulator-private facts.

The memory update is abstracted as $\mathcal{M}_t = \Phi(b,\tau_t,\mathcal{M}_{t-1})$, 
where $\Phi$ is a structured LLM call followed by deterministic validation. Given the public brief, the public transcript, and the previous ledger, it searches for missing requirements across six categories: objectives and trade-offs, decision scope, operational constraints, time boundaries, relationships among entities or decisions, and hard-versus-soft policies. It proposes at most three additions, each carrying a category, a description, and supporting public evidence. Exact duplicates after text normalization are dropped; the remaining additions receive stable identifiers and the status \textsc{Open}. This stage only records gaps: it neither ranks them nor asks the user a question.

The currently bindable gaps are
\begin{equation}
\mathcal{O}_t = \{\ell \in \mathcal{M}_t :
\operatorname{status}(\ell)=\textsc{Open}\}.
\end{equation}
Thus, $\mathcal{O}_t$ contains registered gaps that have not yet been asked. An empty $\mathcal{O}_t$ does not establish that the formulation is complete.

\subsection{Stage 2: Gap-Guided Action Search}

At each turn, Stage 2 runs $\Psi$ on the public brief, the public transcript, and $\mathcal{O}_t$.
$\Psi$ returns a stopping decision $d_t$ together with three candidate clarification states $C_t$.
Each state records the confirmed goal, decision scope, constraints and rules, known entities and inputs, and unresolved business assumptions.
If $d_t$ is \textsc{Ask}, $\Psi$ pairs each state with one candidate question, yielding three corresponding questions $Q_t$; under the Choice protocol, each question offers exactly three answer options.

Whenever $\mathcal{O}_t \neq \emptyset$, every candidate must name the identifier of the open entry it targets. If no open entry remains, the candidates may address another clarification need without a ledger binding. The runner validates the output schema and any required gap binding, regenerating invalid output up to three times.

When $\Psi$ decides to ask, a selector compares the three state--question pairs against the public history and open ledger:
\begin{equation}
a_t =
\operatorname{Select}
(\tau_t,\mathcal{C}_t,\mathcal{Q}_t,\mathcal{O}_t).
\end{equation}
The selector is instructed to consider the candidates in the following order: potential changes to objectives or trade-offs, business constraints and rules, decision scope, and the risk of a silent assumption. It downweights candidates concerned only with mathematical detail. Only after the chosen question enters the public transcript does its bound entry move from \textsc{Open} to \textsc{Asked}, marking that the gap was queried, not resolved. In the choice-based instantiation of \methodname{}, each selected question is presented with candidate answer options.

\subsection{Algorithm}

Algorithm~\ref{alg:interopt} summarizes the full procedure.

\begin{algorithm}[H]
\caption{\methodname{} Clarification}
\label{alg:interopt}
\begin{algorithmic}[1]
\REQUIRE Public brief $b$, protocol $\pi$, maximum turns $T_{\max}$
\STATE $\tau_1 \leftarrow \emptyset$, $\mathcal{M}_0 \leftarrow \emptyset$
\FOR{$t=1$ to $T_{\max}$}
  \STATE $\mathcal{M}_t \leftarrow \Phi(b,\tau_t,\mathcal{M}_{t-1})$
  \STATE $\mathcal{O}_t \leftarrow \{\ell \in \mathcal{M}_t : \operatorname{status}(\ell)=\textsc{Open}\}$
  \STATE $(d_t,\mathcal{C}_t,\mathcal{Q}_t) \leftarrow \Psi(b,\tau_t,\mathcal{O}_t,\pi)$
  \IF{$d_t = \omega$}
    \RETURN $\tau_t$ and $\omega$
  \ENDIF
\STATE Constrain $\mathcal{Q}_t$ to $\mathcal{O}_t$ when $\mathcal{O}_t\neq\emptyset$
  \STATE $a_t \leftarrow \operatorname{Select}(\tau_t,\mathcal{C}_t,\mathcal{Q}_t,\mathcal{O}_t)$
  \STATE Execute public action $a_t$ under protocol $\pi$
  \STATE Mark the gap bound to $a_t$ as \textsc{Asked}, if any
  \STATE Receive fact-bounded simulated-user answer $u_t$
  \STATE $\tau_{t+1} \leftarrow \tau_t \cup \{(a_t,u_t)\}$
\ENDFOR
\RETURN $\tau_{T_{\max}+1}$ with turn-limit termination
\end{algorithmic}
\end{algorithm}
\FloatBarrier
\section{Experiments and Results}
\label{sec:experiments}

\subsection{Experimental Setup}
\label{sec:experimental-setup}

All main experiments use the 100 cases in \benchmarkname{} with five independent runs per case unless otherwise stated.
We treat the case as the statistical unit by averaging repeated runs within each case before aggregating across cases.
Runs that fail the post-hoc protocol audit remain in the headline metric denominators and are retained for diagnostic analysis.
We report All-Slot Exact, Core Exact, Silent/run, Avg Turns, and Avg Q.
Core Exact follows Section~\ref{sec:benchmark} and includes only cases with at least one $P0/P1$ slot (94 of 100 cases), while All-Slot Exact includes all cases; partial recovery is not counted as exact recovery.
Silent/run counts unconfirmed hidden-slot assumptions, and Avg Turns/Avg Q measure interaction length at the turn and atomic-question levels.
Unless otherwise specified, the method-comparison and ablation experiments use DeepSeek V4 Pro as the tested model, $T_{\max}=20$, agent temperature 0.2, and simulator/judge/selector temperature 0.0; for \methodname{}, Stage 1 adds at most three new gaps per turn and Stage 2 instantiates three candidate states and questions on each asking turn.

For methods with internal planning roles, all agent-facing decisions--gap discovery, ask/ready decision, candidate generation, and candidate selection--are instantiated with the tested agent model. The simulated user, protocol detector, and post-hoc judge are fixed evaluation components and are not allowed to expose hidden slots or recovery labels to the tested agent. We therefore evaluate a prompt-level interaction policy under a fixed harness.

\subsection{Clarification Protocols}
\label{sec:clarification-protocols}

\subsubsection{Choice-Based Setting}
\label{sec:choice-experiments}

We compare four methods in the Choice block, each evaluated with \(K=5\) runs per case.
\textbf{MC} asks questions with three generated options A--C; when none matches, the simulated user must choose the closest option.
\textbf{MC-D} uses the same agent but adds a fixed option D that permits a free-form correction.
MC and MC-D are controlled protocols designed for this study. They instantiate option-based clarification for pre-formulation clarification, drawing on prior QA work with user-selectable alternatives \citep{lee2023cambignq}.
MC-D serves as the base interaction protocol for the Choice-based variants built on it.
\textbf{ReadyGate} augments MC-D with an independent stopping reviewer that accepts or rejects \readytoken{} without seeing hidden slots.
\textbf{\methodname{}} is the complete two-stage method: Stage 1 maintains a persistent gap memory; Stage 2 anchors candidate questions to open gaps when available, while leaving the ask-or-stop decision to the agent.

When the agent asks a question, it generates options A--C, and the MC-D runner appends the fixed option D. The simulated user selects the best-fitting option among A--C; when those options do not represent the relevant private fact, the response is submitted in free form through D. Internally, the simulator records a match-status label (exact, acceptable, no-match, or undetermined) solely for diagnostics and does not reveal it to the agent or judge.
The protocol monitor passively validates each action without intervening; it does not reject multi-question turns or trigger regeneration.
All Choice methods share the same cases, repetitions, model, simulator, judge, and monitor.

\subsubsection{Open-Ended Setting}

The open/free-form comparison includes exactly the methods reported in the Open block of Table~\ref{tab:resource_results}.
All methods use the same FreeQA interface: each action is a natural-language question or \readytoken{}, and all transcripts are scored by the same judge.
\textbf{FreeQA} is our base open protocol without an external planner or stopping reviewer.
\textbf{ReadyGate} adds a readiness reviewer that checks unresolved assumptions when the agent attempts to stop.
\textbf{\methodname{}} is our two-stage formulation-gap-guided framework. 

The remaining Open baselines are task-aligned adaptations of the original systems. Full end-to-end reproduction would require components beyond pre-formulation clarification. The GATE adapter retains the informative-question principle and uses \readytoken{} as the stopping action. The ORPilot adapter retains its interview stage for eliciting objectives, decisions, constraints, parameters, and indices; downstream data collection, code generation, solver execution, and reporting are outside our evaluation. The Let's T-P-P adapter retains its domain-prompted, process-aware conversational policy while excluding the task-specific SFUSD environment, tool calls, solver loop, and hidden-utility evaluation. We cite the corresponding original systems as \citep{tamkin2023gate,xie2026orpilot,drossman2026conversation}.

For controlled comparison, all Open adapters receive the same public brief and dialogue history, use the same FreeQA action space and $T_{\max}=20$, and are evaluated on the same 100 cases with $K=5$ runs, the same tested model, simulated user, and post-hoc judge. Neither the adapters nor their stopping actions receive hidden slots or recovery rubrics.

In the open/free-form setting, the agent asks natural-language questions or emits \readytoken{}.
The simulated user answers only from the private case facts and dialogue history, without volunteering unasked hidden facts.
The monitor records invalid actions for diagnostics.
Atomic questions are counted separately from turns because one utterance may contain multiple questions.

\subsection{Performance across Tested Models}

\begin{table*}[!ht]
\centering
\caption{Baseline model comparison in \benchmarkname{}. 
}
\label{tab:main_results}
\small
\setlength{\tabcolsep}{5pt}

\begin{tabular}{llrrrrr}
\toprule
Setting & Tested model & All-Slot Exact $\uparrow$ & Core Exact $\uparrow$ &
Silent/run $\downarrow$ & Avg Turns & Avg Q \\
\midrule
Open / FreeQA & DeepSeek V4 Pro & 0.426 & 0.472 & 0.682 & 3.912 & 3.232 \\
Open / FreeQA & GLM-5.1         & 0.454 & 0.491 & 0.624 & 4.484 & 3.880 \\
Open / FreeQA & GPT-5.5         & 0.430 & 0.498 & 0.632 & 3.491 & 2.540 \\
Open / FreeQA & Opus-4.8        & \textbf{0.548} & \textbf{0.583} & \textbf{0.424} & 4.406 & 3.504 \\

\midrule
Choice / MC-D & DeepSeek V4 Pro & 0.474 & 0.506 & 0.692 & 3.258 & 2.350 \\
Choice / MC-D & GLM-5.1         & 0.400 & 0.434 & 0.712 & 2.898 & 1.970 \\
Choice / MC-D & GPT-5.5         & 0.414 & 0.449 & 0.660 & 3.184 & 2.218 \\
Choice / MC-D & Opus-4.8        & \textbf{0.542} & \textbf{0.583} & \textbf{0.456} & 3.724 & 2.906 \\
\bottomrule
\end{tabular}

\vspace{0.3em}
\end{table*}

Table~\ref{tab:main_results} compares off-the-shelf LLMs under Open/FreeQA and Choice/MC-D. Performance varies by model and protocol; Opus-4.8 is strongest in both settings, but no model exceeds 60\% Core Exact and all retain substantial silent assumptions. These results confirm that OR-Clarify remains challenging.

\subsection{Method Effectiveness and Interaction Cost}

\begin{table*}[!ht]
\centering
\caption{Training-free method comparison and interaction burden in \benchmarkname{}.}
\label{tab:resource_results}
\small
\begin{tabular}{llrrrrr}
\toprule
Setting & Method & All-Slot Exact $\uparrow$ & Core Exact $\uparrow$ & Silent/run $\downarrow$ & Avg Turns & Avg Q \\
\midrule
Open & FreeQA        & 0.426 & 0.472 & 0.682 & 3.912 & 3.232 \\
Open & ReadyGate     & 0.452 & 0.494 & 0.418 & 6.476 & 5.576 \\
Open & ORPilot       & \textbf{0.546} & \textbf{0.583} & \textbf{0.268} & 5.914 & 4.918 \\
Open & Let's T-P-P    & 0.458 & 0.489 & 0.640 & 4.756 & 3.766 \\
Open & GATE          & 0.528 & 0.560 & 0.412 & 5.544 & 4.818 \\
Open & \methodname{} & 0.492 & 0.538 & 0.560 & 5.490 & 4.824 \\
\midrule
Choice & MC            & 0.414 & 0.447 & 0.774 & 3.326 & 2.400 \\
Choice & MC-D          & 0.474 & 0.506 & 0.692 & 3.258 & 2.350 \\
Choice & ReadyGate     & 0.476 & 0.517 & 0.614 & 3.982 & 2.892 \\
Choice & \methodname{} & \textbf{0.638} & \textbf{0.675} & \textbf{0.366} & 10.674 & 10.036 \\
\bottomrule
\end{tabular}
\end{table*}

Table~\ref{tab:resource_results} compares training-free clarification methods under both interfaces. In the open/free-form setting, no method is uniformly dominant: ORPilot attains the strongest recovery, GATE remains competitive, and \methodname{} is close but not best while still leaving more silent assumptions than the strongest open baselines. This shows that formulation-gap-guided clarification does not uniformly dominate strong free-form baselines.

\begin{figure}[!t]
\centering
\includegraphics[width=1\columnwidth]
{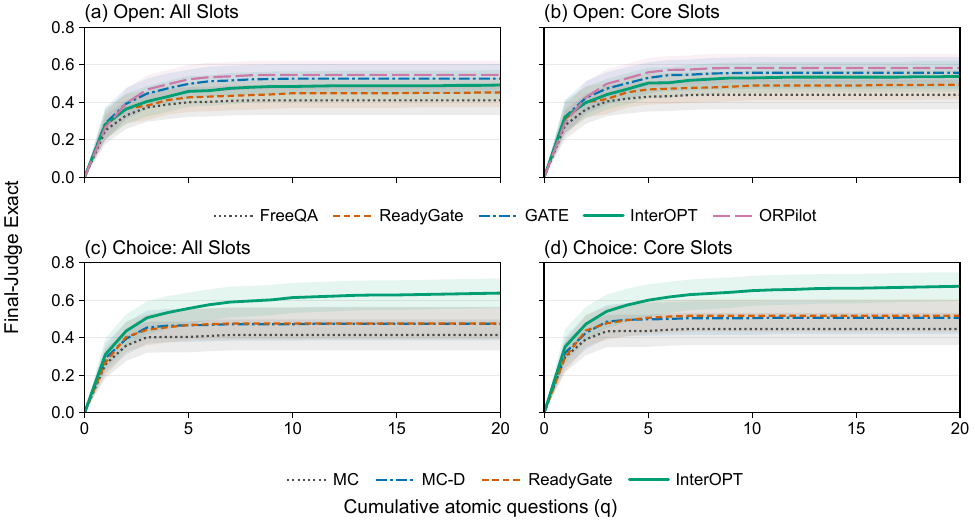}
\caption{Exact restoration versus cumulative atomic questions under Open and
Choice protocols. Top: Open; bottom: Choice. Left: All-Slot Exact; right:
Core Exact. Curves average $K=5$ runs, with 95\% case-level bootstrap confidence
intervals. Methods are compared only within the same response protocol.}
\label{fig:question_efficiency}
\end{figure}

In the Choice setting, \methodname{} obtains the highest recovery among the Choice methods we evaluate, while also using substantially more turns and atomic questions. The overall pattern is therefore a coverage–cost trade-off: structured gap-guided clarification improves recovery most under constrained answer spaces, but stopping and question efficiency remain limiting factors.

Table~\ref{tab:resource_results} reports endpoint performance under each method's own stopping policy. It therefore answers two questions: how much requirement recovery a complete run ultimately obtains, and how many questions the stopping policy spends. Figure~\ref{fig:question_efficiency} provides the complementary budget-matched view. At a cumulative budget of $q$ atomic questions, each point reports the case-averaged exact-recovery rate achieved by that point; methods are compared at the same $q$ and only within the same interaction protocol. The curves average $K=5$ runs within each case, and the shaded 95\% case-level bootstrap intervals reflect variation across cases.

Under Choice, \methodname{} continues gaining exact recovery over later questions after MC-D and ReadyGate largely plateau. Its endpoint advantage should therefore be read together with the overlapping portions of the curves, where recovery is compared under a common question budget. 
The Open panels retain ORPilot and the other free-form baselines under their own shared response protocol.

\subsection{Contributions of the Two Stages}

We evaluate two system-level ablations under both protocols. \emph{w/o Stage 1} disables formulation-gap search and ledger population while retaining Stage-2 candidate states, three candidate questions, and the selector. \emph{w/o Stage 2} retains gap search and the persistent ledger but collapses the architecture into a direct gap-to-question policy, with Stage 1 supplying both the public question and the readiness decision. Both variants use the same cases, repetitions, temperatures, simulator, and judge as full \methodname{}.

\begin{table*}[!ht]
\centering
\caption{System-level ablation of \methodname{} under open/FreeQA and Choice (MC-D) clarification.}
\label{tab:interopt_ablation}
\small
\begin{tabular}{llrrrrr}
\toprule
Clarification & Variant & All-Slot Exact $\uparrow$ & Core Exact $\uparrow$ &
Silent/run $\downarrow$ & Avg Turns & Avg Q \\
\midrule
Open / FreeQA & None          & 0.426 & 0.472 & 0.682 & 3.912 & 3.232 \\
Open / FreeQA & w/o Stage 1   & 0.432 & 0.460 & 0.668 & 4.074 & 3.492 \\
Open / FreeQA & w/o Stage 2   & 0.470 & 0.506 & \textbf{0.254} & 9.514 & 9.502 \\
Open / FreeQA & \methodname{} & \textbf{0.492} & \textbf{0.538} & 0.560 & 5.490 & 4.824 \\
\midrule
Choice (MC-D) & None          & 0.474 & 0.506 & 0.692 & 3.258 & 2.350 \\
Choice (MC-D) & w/o Stage 1   & 0.520 & 0.564 & 0.612 & 4.346 & 3.456 \\
Choice (MC-D) & w/o Stage 2   & 0.572 & 0.615 & 0.462 & 5.700 & 4.860 \\
Choice (MC-D) & \methodname{} & \textbf{0.638} & \textbf{0.675} & \textbf{0.366} & 10.674 & 10.036 \\
\bottomrule
\end{tabular}

\vspace{0.3em}
\end{table*}

Table~\ref{tab:interopt_ablation} reports both ablations. Under the open protocol, removing Stage 1 yields 0.460 Core Exact, whereas removing Stage 2 yields 0.506 but asks nearly twice as many questions as full \methodname{} (9.5 vs.\ 4.8). Within these system-level ablations, Stage 1 has the larger contribution to open-setting recovery, whereas Stage 2 mainly reduces interaction length. Under the Choice protocol, both stages contribute, with the removal of Stage 1 causing the larger reduction (11.1 vs. 6.0 percentage points in Core Exact). This pattern suggests a division of labor between the two stages: gap discovery contributes most directly to coverage, while guided action selection becomes more consequential for recovery when clarification is conducted through the structured Choice interface.

\subsection{Protocol-Specific Diagnostics}

\paragraph{Choice diagnostics.}
In MC-D, all 132 no-match events (among 1,129 audited choice events) invoked D and exposed a public free-form correction.
These descriptive event rates confirm the protocol distinction, but they are not a causal estimate of D because the two agents can generate different question trajectories.

\methodname{} leaves 0.366 silent assumptions per run, compared with 0.692 for MC-D and 0.614 for ReadyGate.
Exact restoration continues to treat partial matches as failures: across 890 judged slots per $K=5$ method, \methodname{} receives 617 \texttt{yes}, 20 \texttt{partial}, and 253 \texttt{no} labels, whereas MC-D receives 461, 16, and 413.
Thus, the headline gain comes primarily from converting unresolved slots into fully restored ones, with little contribution from partial-credit cases. The small number of partial labels further suggests that the improvement is concentrated on explicitly reaching the relevant requirement, with relatively few gains coming from near-miss clarification questions.

\paragraph{Open diagnostics.}
Open-specific diagnostics separate question packaging from stopping errors. Packaging is not the explanation: \methodname{} asks 4.824 atomic questions over 4.502 question turns per run, and only 6.4\% of question turns bundle more than one question. The dominant failure is stopping under residual formulation uncertainty—the agent declares readiness in 98.8\% of runs, yet the stopping audit flags 203 of 500 runs (40.6\%) as premature, which accounts for the 0.560 silent assumptions that remain. Nor does recovery depend on unusually informative answers: multi-slot disclosures occur only 0.116 times per run, so the gap memory accumulates evidence across turns as intended. The open-setting challenge is thus not that the agent fails to ask, but that it fails to know when it has asked enough.
\section{Conclusion}

We introduced \benchmarkname{}, a benchmark for pre-formulation clarification in operations research, and \methodname{}, a two-stage framework that maintains formulation gaps to guide questioning and stopping. \benchmarkname{} evaluates whether agents can recover formulation-critical hidden requirements before modeling, while also auditing interaction cost, stopping behavior, and silent assumptions.

Empirically, \methodname{} is most effective in the Choice setting, where structured gap-guided clarification improves exact recovery over training-free baselines. In the open/free-form setting, it remains competitive but does not uniformly dominate strong baselines, showing that its gains do not transfer uniformly to free-form clarification. Ablations further indicate that both gap tracking and guided action selection contribute to recovery.

These results also expose open challenges. \methodname{} can require more interaction, stopping remains imperfect, and silent assumptions are not eliminated. Future work should expand the benchmark, improve question efficiency, and develop better stopping calibration. Overall, optimization agents should be evaluated not only after they produce a model, but also before modeling: on whether they know what to ask, when to ask, and when the specification is complete enough to model.

\small
\bibliographystyle{plainnat}
\bibliography{refs}

\end{document}